\documentclass{IEEEtran}

\IEEEoverridecommandlockouts                              

\usepackage{graphics} 
\usepackage{epsfig} 
\usepackage{mathptmx} 
\usepackage{times} 
\usepackage{amsmath} 
\usepackage{amssymb}  
\usepackage{float}
\usepackage{color}
\usepackage{algorithm}
\usepackage{algorithmic}

\makeatletter
\def\widebreve#1{\mathop{\vbox{\m@th\ialign{##\crcr\noalign{\kern3\p@}%
      \brevefill\crcr\noalign{\kern3\p@\nointerlineskip}%
      $\hfil\displaystyle{#1}\hfil$\crcr}}}\limits}

\def\brevefill{$\m@th \setbox\z@\hbox{$\braceld$}%
  \bracelu\leaders\vrule \@height\ht\z@ \@depth\z@\hfill\braceru$}
\makeatletter

\usepackage{cite}
\usepackage{amsfonts}
\usepackage{graphicx}
\usepackage{textcomp}

\def\BibTeX{{\rm B\kern-.05em{\sc i\kern-.025em b}\kern-.08em
    T\kern-.1667em\lower.7ex\hbox{E}\kern-.125emX}}

\begin{document}



\title{\LARGE \bf
Optimal Strategies for Disjunctive Sensing and Control
}

\author{Richard L. Sutherland$^{1}$, Ilya V. Kolmanovsky$^{2}$, Anouck R. Girard$^{3}$, \\ Frederick A. Leve$^{4}$, and Christopher D. Petersen$^{5}$
\thanks{$^{1}$PhD Candidate, Department of Aerospace Engineering,
        University of Michigan, Ann Arbor, MI 48109, USA
        {\tt\small rlsu@umich.edu}}%
\thanks{$^{2}$Professor, Department of Aerospace Engineering, University of Michigan,
        Ann Arbor, MI 48109, USA
        {\tt\small ilya@umich.edu}}%
\thanks{$^{3}$Associate Professor, Department of Aerospace Engineering, University of Michigan,
        Ann Arbor, MI 48109, USA
        {\tt\small anouck@umich.edu}}%
\thanks{$^{4}$Program Officer, Air Force Office of Scientific Research,
        Arlington, VA 22203, USA
        {\tt\small frederick.leve@us.af.mil}}%
\thanks{$^{5}$Research Engineer, Space Vehicles Directorate
        {\tt\small AFRL/VSSVOrgMailbox@us.af.mil}}%
}


\maketitle
\thispagestyle{empty}
\pagestyle{empty}

\begin{abstract}
A disjunctive sensing and actuation problem is considered in which the actuators and sensors are prevented from operating together over any given time step. This problem is motivated by practical applications in the area of spacecraft control. Assuming a linear system model with stochastic process disturbance and measurement noise, a procedure to construct a periodic sequence that ensures bounded states and estimation error covariance is described along with supporting analysis results. The procedure is also extended to ensure eventual satisfaction of probabilistic chance constraints on the state. The proposed scheme  demonstrates good performance in simulations for spacecraft relative motion control.



\end{abstract}


\section{Introduction}
\label{sec:intro}


A common assumption in control theory is that sensing and actuation can be performed simultaneously. In this paper, we consider the case of \textit{disjunctive} sensing and actuation, in which, at any given time step, either a sensor or an actuator can be operated but not both. Thus, a switching policy between sensing and actuation needs to be determined that achieves the specified mission objectives.  

The motivation for considering this class of problems comes from spacecraft control applications. Specifically, magnetic fields generated by magnetic actuators may interfere with the magnetic sensors used for attitude sensing, see e.g., \cite{Mazzini}. In larger satellites, magnetometers may be placed on a boom to reduce the electromagnetic interference from the magnetic actuators and other equipment onboard the spacecraft. Such a solution is not feasible for smaller and cheaper cubesats, where the magnetic actuators must be deactivated to record an accurate attitude reading. Other situations in which simultaneous sensing and actuation are not possible include vision-based rendezvous, docking, and relative motion maneuvering \cite{Fehse}, in which the plume or vibration from spacecraft thrusters may interfere with cameras or other sensitive navigation equipment.

For this problem, we consider the eventual enforcement of probabilistic chance state constraints of the form
\begin{equation}
    \forall k \geq k^*, \ \text{Prob}(\{ x_k \in X \}) \geq 1 - \delta,
\label{eqn:cnc_con}
\end{equation}
where $X$ is the set prescribed by the constraints of the problem, $0 \leq \delta < 1$, and $k^* \in \mathbb{Z}_+$ is sufficiently large. In spacecraft applications, $X$ may represent a region in the state space in which scientific measurements can be reliably taken by the onboard instrumentation, see, for instance, the case study in \cite{SFM} and \cite{TCST_arxiv}. In this case study, the state constraints do not need to be satisfied during initial transients but must be satisfied eventually to enable the equipment to function.

In this paper, we treat a disjunctive sensing and actuation problem for systems that can be represented by discrete-time linear models with stochastic process disturbance and measurement noise inputs. A procedure to construct a periodic switching sequence between sensing and actuation is described and closed-loop boundedness and convergence properties when such a sequence is applied are analyzed.


The given problem falls within the general class of switched system stabilization problems in which a part of the dynamics represents the propagation of the estimated state and of the estimation error covariance matrix. Stability of switched systems has been studied extensively, see e.g., \cite{Lieberzon} and \cite{Sun}.  We note that if the system is open-loop unstable then, to guarantee either boundedness of states or boundedness of state estimates based on dwell time conditions, sufficient dwell time in each mode (sensing or actuation) is necessary. In addition, sensing and actuation intervals need to be suitably interlaced to achieve both goals simultaneously.

Techniques developed for stability analysis and control of discrete-time periodic systems \cite{Bittanti,Sadkane,Zhou} are also relevant given our search for a periodic switching sequence. In event triggered and self-triggered control \cite{Heemels} and in sensor networks, sensor scheduling and sensor tasking \cite{Liu,Sunberg}, sensors may be deactivated when confidence in the estimated states is high; however, the situation in which the sensors and the actuators cannot be used simultaneously does not appear to be treated.

The main contributions of the paper are the formulation of a disjunctive sensing and actuation problem based on linear discrete-time system models and a procedure to construct a periodic switching sequence. The constructed sequence ensures the convergence of the state mean, the boundedness of the estimation error covariance matrix, and the eventual satisfaction of the probabilistic chance state constraints. Simulation examples are reported that confirm good behavior of the proposed scheme for a spacecraft relative motion control example.

The paper is organized as follows. In Section~\ref{sec:dsc}, a disjunctive sensing and actuation problem is formally stated. In Section~\ref{sec:swtch_seq}, properties of the closed-loop system operated with periodic switching sequences are analyzed and admissible sequences that lead to desired boundedness and convergence properties are characterized. In Section~\ref{sec:chancecon}, the eventual satisfaction of the probabilistic chance state constraints is described. The selection of an optimal admissible periodic switching sequence is discussed in Section~\ref{sec:opt}. A simulation case study that addresses the relative motion control problem is presented in Section~\ref{sec:opt_cost_ex}. Finally, concluding remarks are made in Section~\ref{sec:concl}.


The notations used are standard: $\mathbb{E}$ denotes the expectation, $\mathbb{Z}_{\geq 0}$ denotes the set of nonnegative integers, $\mathbb{Z}_{[0,N]}$ denotes the set of nonnegative integers between $0$ and $N$, Tr$(P)$ denotes the trace of a matrix $P$, $\| Q \|$ denotes a norm of a matrix $Q$, and $\rho(R)$ denotes the spectral radius of a matrix $R$. For two symmetric matrices $(P, \ Q)$, $P \succ Q$ implies that $P-Q$ is positive definite and $P \succeq Q$ implies that $P-Q$ is positive semi-definite. The matrix $\mathbf{I}_n$ is the $n \times n$ identity matrix and the matrix $\mathbf{0}_{m \times n}$ is the $m \times n$ zero matrix; we drop subscripts in these when they are clear from the context.

A conference version appeared in the 2018 ACC proceedings \cite{ACC2018}. As compared to the conference paper, this version contains proofs and additional details not found in the conference paper such as an expanded treatment of chance constraints.

\section{Disjunctive Sensing and Actuation for Linear Systems}
\label{sec:dsc}

Consider a system represented by a discrete-time linear model with stochastic state disturbance and measurement noise inputs, given by
\begin{equation}
\begin{split}
    x_{k+1} &= A x_k + B u_k + w_k, \\
    y_k &= C x_k + \nu_k,
\label{eqn:state}
\end{split}
\end{equation}
where $x_k$ is a vector state and $u_k$ is a vector control. The variables $w_k$ and $\nu_k$ are, respectively, the state disturbance and the measurement noise inputs, each assumed to be a sequence of zero-mean, independent (and jointly independent) identically distributed (i.i.d.) random variables with $\mathbb{E}[w_k w_k^{\rm T}] = \Sigma_w = \Sigma_w^{\rm T}$, $\mathbb{E}[\nu_k \nu_k^{\rm T}] = \Sigma_\nu = \Sigma_\nu^{\rm T}$.

In this paper, a fixed gain state feedback is considered,
\begin{equation}
    u_k = u_T + K (\hat{x}_k - x_T),
\label{eqn:feedback}
\end{equation}
where $\hat{x}_k$ is the state estimate generated by a fixed gain observer of the form
\begin{equation}
    \hat{x}_{k+1} = A \hat{x}_k + \eta_k B u_k + (1 - \eta_k) L (y_k - C \hat{x}_k).
\label{eqn:observer}
\end{equation}
The binary variable $\eta_k \in\{0,1\}$ represents the system operating mode. When $\eta_k = 1$, the control is applied to the system.  When $\eta_k = 0$, the control is deactivated ($u_k = 0$) and the sensed output, $y_k$, is obtained. The target equilibrium is denoted by $x_T$, and it is assumed that the feedforward control input, $u_T$, supports it in steady-state in the absence of $w_k$ and with $\eta_k = 1$, i.e., $x_T = Ax_T + Bu_T$.

In a typical control design process, due to the Separation Principle, the gains $K$ and $L$ would be determined without consideration of the interference between actuation and sensing (possibly by different engineers) and then a coordination mechanism introduced by specifying $\eta_k, \ k = 0, 1, \dots$. In this setting, offline-generated $N$-periodic switching sequences, $\{\eta_k\}$, with $\eta_k = \eta_{k+N}$ for all $k\in \mathbb{Z}_{\geq 0}$, $N \in \mathbb{Z}_{>0}$, are of particular interest, so that their repeated application leads to the attainment of the control objectives, including eventual satisfaction of state constraints.
This approach, based on the application of the offline generated periodic sequence, has low computational footprint and is appealing in view of limited computing power and restrictive electrical power consumption budgets onboard of small spacecraft.



\section{Admissible Periodic Switching Sequences}
\label{sec:swtch_seq}

Suppose that $\{ \eta_k \}$ is fixed and define $\bar{A}_k = A + \eta_k B K, \ ~\tilde{A}_k = A + (1-\eta_k)LC$. The evolution of the state, $x_k$, and of the state estimation error, $e_k = x_k - \hat{x}_k$, are driven by
\begin{equation}
\begin{split}
        x_{k+1} &= A x_k + \eta_k B K \hat{x}_k + w_k \\
                &= \bar{A}_k x_k - \eta_k B e_k + \eta_k B u_T + w_k,
    \end{split}
\label{eqn:state_dyn1}
\end{equation}
and
\begin{equation}
    \begin{split}
        e_{k+1} &= \left[ A + (1 - \eta_k) L C \right] e_k + w_k + (1 - \eta_k)L \nu_k \\
                &= \tilde{A}_k e_k + w_k + (1 - \eta_k)L \nu_k.
    \end{split}
\label{eqn:est_err}
\end{equation}
In simulations, we assume that $x_T = 0$ (and so $u_T = 0$). Based on (\ref{eqn:est_err}) and assumed independence and zero mean properties of the stochastic disturbance, $w_k$, and measurement noise, $\nu_k$, the error covariance matrix $P_{k} = \mathbb{E}[e_k e_k^\textrm{T}]$ satisfies
\begin{equation}\label{eqn:Pn}
    P_{k+1} = \tilde{A}_k P_{k} \tilde{A}_k^\textrm{T} + R_k,
\end{equation}
\begin{equation}\label{eqn:Rn}
R_k =
    \left[\begin{array}{cc}(1-\eta_k)L & {\bf I} \end{array}\right] \left[\begin{array}{cc} \Sigma_\nu & {\bf 0} \\ {\bf 0} & \Sigma_w \end{array} \right] \left[\begin{array}{c} (1-\eta_k)L^{\rm T} \\ {\bf I} \end{array} \right].
\end{equation}

\textbf{Definition 1:} Let $A$, $B$ be defined as in (\ref{eqn:state}), $K$, $L$ be defined as in (\ref{eqn:feedback}), (\ref{eqn:observer}), and $\bar{A}_k$, $\tilde{A}_k$ be defined as above, with none of $\bar{A}_k$, $\tilde{A}_k$ nilpotent. An $N$-periodic sequence of binary integers $\{ \eta_0, \eta_1, \cdots, \eta_{N-1} \}$, where, for all $k \in \mathbb{Z}_{\geq 0}$, $\eta_k \in \{0,1\}$ and $\eta_{k+N}=\eta_k$, is called \textit{admissible} if the following contractivity conditions hold:
\begin{equation}
    \rho (\bar{A}_{N-1} \bar{A}_{N-2} \cdots \bar{A}_0) = \bar{q}_{A} < 1,
\label{eqn:def_1a}
\end{equation}
\begin{equation}
    \rho (\tilde{A}_{N-1} \tilde{A}_{N-2} \cdots \tilde{A}_0) = \tilde{q}_{A} < 1,
\label{eqn:def_1b}
\end{equation}
where $\rho(\cdot)$ denotes the spectral radius operator. A non-periodic sequence of binary integers $\{ \eta_0, \eta_1, \dots \}$ is admissible if:
\begin{equation}
    \lim_{k \to \infty} \rho (\bar{A}_{k} \bar{A}_{k-1} \cdots \bar{A}_0) = 0,
\label{eqn:def_1c}
\end{equation}
\begin{equation}
    \lim_{k \to \infty} \rho (\tilde{A}_{k} \tilde{A}_{k-1} \cdots \tilde{A}_0) = 0.
\label{eqn:def_1d}
\end{equation}
The constraint against nilpotent matrices ensures that the limits in (\ref{eqn:def_1c}) and (\ref{eqn:def_1d}) approach zero rather than ``jump'' to zero. These definitions are consistent with the properties of discrete-time state transition matrices found in, for example, Chen \cite{Chen}.

\textbf{Lemma 1:} If an admissible sequence exists, then a periodic admissible sequence exists.

\textit{Proof:} Let $s = \{ \eta_0, \eta_1, \dots, \eta_k, \dots \}$ be an admissible sequence. If $s$ is periodic, then done. Thus, assume $s$ is non-periodic. By definition, $\lim_{k \to \infty} \rho (\bar{A}_{k} \bar{A}_{k-1} \cdots \bar{A}_0) = 0,$
thus, for each $\bar{\epsilon} > 0$, there exists $\bar{N} \in \mathbb{N}$ such that, for all $k > \bar{N}, ~\rho (\bar{A}_{k} \bar{A}_{k-1} \cdots \bar{A}_0) < \bar{\epsilon},$
and, by similar argument, for each $\tilde{\epsilon} > 0$ there exists $\tilde{N}$ such that, for all $k > \tilde{N}, ~\rho (\tilde{A}_{k} \tilde{A}_{k-1} \cdots \tilde{A}_0) < \tilde{\epsilon}$. Choose $k$ such that max$\{\bar{\epsilon}, \tilde{\epsilon}\} < 1$. Then, $\rho (\bar{A}_{k} \bar{A}_{k-1} \cdots \bar{A}_0) < 1$, and $\rho (\tilde{A}_{k} \tilde{A}_{k-1} \cdots \tilde{A}_0) < 1$, therefore $s_k = \{ \eta_0, \eta_1, \dots, \eta_k\}$ is admissible by construction and, when applied repeatedly, is periodic with period $k+1$. $\blacksquare$


\subsection{Limits of the Mean and Error Covariance Matrix Sequences}
\label{sec:covar}
Note first that $\mathbb{E}[\nu_k] = 0$, $\mathbb{E}[w_k] = 0$ and hence (\ref{eqn:state_dyn1}), (\ref{eqn:est_err}) imply that the state and estimation error mean, $\mu_{x,k} = \mathbb{E}[x_k]$ and $\mu_{e,k}=\mathbb{E}[e_k]$, respectively, satisfy
\begin{equation}
    \mu_{x,k+1} = \bar{A}_k \mu_{x,k} - \eta_k B \mu_{e,k} + \eta_k B u_T,
\label{eqn:mean_x}
\end{equation}
\begin{equation}
    \mu_{e,k+1} = \tilde{A}_k \mu_{e,k},
\label{eqn:mean_e}
\end{equation}
where $\bar{A}_k$ and $\tilde{A}_k$ are periodic with the same period, $N$, as $\eta_k$.

{\bf Proposition~1:} Suppose that (\ref{eqn:def_1a}) and (\ref{eqn:def_1b}) hold. Then, as $k \to \infty$, $\mu_{e,k} \to 0$ and $\mu_{x,k} \to {\mu}_{x,k}^{\rm s}$ exponentially, where $\{ {\mu}_{x,k}^s \}$ is the unique $N$-periodic solution of (\ref{eqn:mean_x}) with $\mu_{e,k} \equiv 0$.  Furthermore, if $u_T=0$ then ${\mu}_{x,k}^s=0$. 

\textit{Sketch of the proof:}  The proof follows from Proposition 4.5 in \cite{Halanay} by noting that the characteristic multipliers (eigenvalues of $N$-step state transition matrix) of the combined time-periodic system (\ref{eqn:mean_x})-(\ref{eqn:mean_e}), which is upper triangular, are inside the unit disk if (\ref{eqn:def_1a}) and (\ref{eqn:def_1b}) hold. $\blacksquare$

The next result summarizes the properties of the error covariance matrix sequence.

{\bf Proposition~2:}  Suppose that (\ref{eqn:def_1b}) holds. Then the error covariance matrix, $P_k$, is bounded and, as $k \to \infty$, converges to the unique $N$-periodic solution of (\ref{eqn:Pn}), $\{ P_k^{\rm s} \}$, with $P_{k+N}^{\rm s}=P_k^{\rm s}$. In addition, for any $n \in \mathbb{Z}_{\geq 0}$,
\begin{equation}
    \left( \|P_{N(n+1)}\|-\frac{\gamma}{1-\tilde{q}_A^2} \right) \leq 
    \tilde{q}_A^2 \left( \|P_{Nn}\| - \frac{\gamma}{1-\tilde{q}_A^2} \right),
\label{eqn:Pcontract}
\end{equation}
where 
\[
    \gamma \geq \|\tilde{A}_{N-1}\cdots \tilde{A}_1 R_0 \tilde{A}_1^{\rm T} \cdots \tilde{A}_{N-1}^{\rm T} + \cdots + \tilde{A}_{N-1} R_{N-2} \tilde{A}_{N-1}^{\rm T} + R_{N-1} \|.
\]

\textit{Sketch of the proof:}  The proof of the error covariance matrix convergence follows by applying similar arguments in discrete-time as the ones on p. 58 of \cite{Bittanti} for the continuous-time case. The bound (\ref{eqn:Pcontract}) follows by expressing $P_N$ in terms of $P_0$ and $R_0,\cdots R_{N-1}$ based on (\ref{eqn:Pn}) and applying the triangular inequality. $\blacksquare$


{\bf Remark 1:} The steady-state periodic solution, $P_{k}^s$, of (\ref{eqn:Pn}) can be computed by solving the conventional discrete-time Lyapunov equation for the evolution of the {\em lifted} system error covariance matrix, i.e., of $P^l = \text{diag} \{ P_0^{s}, \cdots, P_{N-1}^s \}$, which is directly obtained from (\ref{eqn:Pn}).

{\bf Remark 2:} Note that (\ref{eqn:Pcontract}) implies that $\limsup_{k \to \infty} \|P_{Nk}\| \leq \gamma/(1-\tilde{q}_A^2)$.

{\bf Remark 3:} The results in Propositions 1 and 2 generalize to non-constant $N$-periodic feedback and observer gains, i.e., $K$ and $L$ are replaced by $K_k$, and $L_k$, where $K_{k+N} = K_k, ~L_{k+N} = L_k$ for all $k \in \mathbb{Z}_{\geq 0}$, under the same conditions (\ref{eqn:def_1a}) and (\ref{eqn:def_1b}). However, analysis results benefit from both $\bar{A}_k$ and $\tilde{{A}}_k$ having only two possible values each, which is the case when the gains are constant.

\subsection{Dwell-Time Conditions and their Implications}

We can take advantage of the dwell time conditions for stability analysis of hybrid systems to develop simpler sufficient conditions that can inform procedures for faster determination of admissible switching sequences. The discussion of the dwell time conditions follows Theorem 4.1 in \cite{Yu} and its proof.

Let $\bar{\Omega}_0 = A$ and $\bar{\Omega}_1 = A + BK$, and consider the condition (\ref{eqn:def_1a}). By Gelfand's theorem \cite{Gelfand}, $\lim_{k \rightarrow \infty} \|\Omega^k\|^{\frac{1}{k}} = \rho(\Omega)$, for any matrix $\Omega$ and norm $\|\cdot\|$. This implies that there exist constants $c_0$, $c_1$, that do not depend on $k$, such that for any $k \geq 1$,
\begin{equation}
    \|\bar{\Omega}_0^k\|^{\frac{1}{k}} \leq c_0 \rho(\bar{\Omega}_0),\quad 
    \|\bar{\Omega}_1^k\|^{\frac{1}{k}} \leq c_1 \rho(\bar{\Omega}_1).
\end{equation}
As the tail of $\|A^k\|^{\frac{1}{k}}$ is strictly non-increasing for any consistent matrix norm, there exists a finite $k^*$ such that $\|A^{k^*}\|^{\frac{1}{k^*}} \geq \|A^k\|^{\frac{1}{k}}$ for all $k \geq 1$. Then, 
\[
c_i = \frac{\|\bar{\Omega}^{k^*}_i\|^{\frac{1}{k^*}}}{\rho(\bar{\Omega}_i)}.
\]
Note that $c_0 \geq 1$ and $c_1 \geq 1$ since $\|A\| \geq \rho(A)$ for any $A$. 
Let $c = \max \{c_0, c_1\}$, $n_0<N$ be the total time spent in the mode $\eta = 0$, $n_1 = N-n_0$ be the total time spent in the mode $\eta = 1$, and $n_s$ be the number of mode ``blocks'' in the sequence, equivalent to the number of switches plus one. Then, $\|\bar{A}_{{N-1}} \bar{A}_{N-2} \cdots \bar{A}_{0} \| \leq c^{n_s} \rho(\bar{\Omega}_0)^{n_0} \rho(\bar{\Omega}_1)^{n_1}$,
and (\ref{eqn:def_1a}) holds if
\begin{equation}\label{eqn:ncc2}
    c^{n_s} \rho(\bar{\Omega}_0)^{n_0} \rho(\bar{\Omega}_1)^{n_1} \leq \bar{q}_{{A}} < 1.
\end{equation}
A frequent situation is that $\bar{\Omega}_0$ (no actuation) is unstable and $\rho(\bar{\Omega}_0)>1$, while $\bar{\Omega}_1$ is stable and $\rho(\bar{\Omega}_1)<1$. 
Then (\ref{eqn:ncc2}) dictates that there must be sufficient time spent in the actuation mode and, furthermore, there must be sufficient dwell time (not too many switches) so that $n_s$ and $c^{n_s}$ are small.

Taking the logarithm of the left hand side of (\ref{eqn:ncc2}), it follows that (\ref{eqn:def_1a}) holds if
\begin{equation}
    n_s \log c + n_0 \log \bar{\rho}_0 + n_1  \log \bar{\rho}_1 < 0,
\label{eqn:ncc2restated}
\end{equation}
where $\bar{\rho}_0=\rho(\bar{\Omega}_0)$, $\bar{\rho}_1=\rho(\bar{\Omega}_1)$. When finding admissible sequences, it is therefore possible to first restrict the search to sequences for which $n_0$, $n_1$ and $n_s$ satisfy the condition (\ref{eqn:ncc2restated}). Note that $\bar{q}_{{A}}$ in (\ref{eqn:ncc2}) is an estimate of the rate of convergence of the state. Hence, ensuring that the left hand side of (\ref{eqn:ncc2restated}) is as negative as possible promotes increasing the convergence rate to $x_T$; this may, however, increase the estimation error.

Similar analysis can be applied in the case of (\ref{eqn:def_1b}). Let $\tilde{\rho}_0=\rho(\tilde{\Omega}_0)$, $\tilde{\rho}_1=\rho(\tilde{\Omega}_1)$ where $\tilde{\Omega}_0 = A + LC$ and $\tilde{\Omega}_1 = A$. Then, $\|\tilde{A}_{N-1} \cdots \tilde{A}_0\| \leq c^{n_s} \rho(\tilde{\Omega}_0)^{n_1} \rho(\tilde{\Omega}_1)^{n_0} \leq \tilde{q}_A < 1$. 
By taking the logarithm of this expression, we can see that (\ref{eqn:def_1b}) holds if
\begin{equation}
    n_s \log c + n_1 \log \tilde{\rho}_0 + n_0  \log \tilde{\rho}_1 < 0 .
\label{eqn:ncc2d2xrestated}
\end{equation}
The condition (\ref{eqn:ncc2d2xrestated}) complements (\ref{eqn:ncc2restated}) and can facilitate the initial fast search for admissible sequences.

\textbf{Remark 4:} Conditions (\ref{eqn:ncc2restated}) and (\ref{eqn:ncc2d2xrestated}) together are sufficient, but not necessary, to also satisfy conditions (\ref{eqn:def_1a}) and (\ref{eqn:def_1b}).

\subsection{Reducible and Irreducible Sequences}
\label{sec:red_seq}

The search for admissible sequences can be made faster by discarding sequences that replicate a known inadmissible subsequence.

\textbf{Definition 2:} A sequence $\{s_N\}$ of length $N \in \mathbb{Z}_{>0}$ is called \textit{reducible} if there exists a subsequence $\{s_k\}$ of length $k \in \mathbb{Z}_{> 0}$ such that $k < N$, $N$ is a multiple of $k$, and $\{s_N\} = \{s_k\} \oplus \{s_k\} \oplus \cdots \oplus \{s_k\},~N/k$ times, where $\oplus$ is used to denote sequence concatenation. A sequence which is not reducible is called \textit{irreducible}.

\textbf{Proposition 3:} Every (non-empty) sequence $\{s_N\}$ contains a unique irreducible subsequence $\{s_n\}$.

\noindent

\textit{Proof:} If $\{s_N\}$ is irreducible, then we are done. Thus, assume $\{s_N\}$ to be reducible. There exists then $k_1 \in \mathbb{N}$ such that $k_1 < N$, $N$ is a multiple of $k_1$, and $\{s_N\} = \{s_{k_1}\} \oplus \hdots \oplus \{s_{k_1}\}$, concatenated $N / k_1$ times. If $\{s_{k_1}\}$ is irreducible, then done, as no shorter subsequences exist and no subsequence longer than $\{s_{k_1}\}$ can be irreducible. If $\{s_{k_1}\}$ is reducible, then there exists $k_2 \in \mathbb{N}$ such that $k_2 < k_1$, $k_1$ is a multiple of $k_2$, and $\{s_{k_1}\} = \{s_{k_2}\} \oplus \hdots \oplus \{s_{k_2}\}$, concatenated $k_1 / k_2$ times. The same argument for $\{s_{k_1}\}$ above now repeats for $\{s_{k_2}\}$. The sequence $\{k_1$, $k_2$, $... \}$ eventually terminates in some $k_n$ as each $k_i$ is strictly smaller than the previous, and reaches an absolute minimum value $k_n \geq 1$ after a finite number of iterations. The process eventually yields an irreducible subsequence $\{s_{k_n}\}$ of $\{s_N\}$ that is of minimum length (no irreducible subsequences of smaller length exist). As this subsequence consists of the first $k_n$ elements of $\{s_N\}$, it is also unique, as only one such subsequence is possible for a given $n$. $\blacksquare$

In practice, we need only focus on these irreducible sub-sequences. This is summarized by the following:

\textbf{Proposition 4:} A binary sequence $\{s_N\}$ is admissible if and only if the associated irreducible subsequence $\{s_n\}$ is admissible.


Any reducible admissible sequence can be formed by propagating an irreducible admissible sequence forward in time. Thus, when investigating sequences of length $N$ for admissibility, all sequences for which we have already evaluated the associated irreducible sub-sequence can be discarded.

To check a sequence $\{ s_N \}$ for reducibility, we employ the following algorithm:

1. Let $D = \{ d : \ d \mid N$ and $d < N \}$, i.e., the set of proper divisors of $N$. 

2. If there exists $d \in D$ such that, for every $1 \leq n \leq N - d$, the sequence satisfies $s_n = s_{n+d}$, where $s_n \in \{s_N\}$, then the sequence is reducible.

3. Otherwise, the sequence is irreducible.



\textbf{Remark 5:} If $N$ is prime, then $D = \{ 1 \}$ and the only reducible sequences are the sequence of zeros and the sequence of ones. Under our initial assumption that both actuation and sensing actions are required, these two sequences can immediately be discarded as inadmissible.

\section{Chance constraints}
\label{sec:chancecon}

Consider now the chance constraint (\ref{eqn:cnc_con}). Define $z = \left[ \begin{array}{cc} x^{\rm T} & e^{\rm T} \end{array} \right]^{\rm T}$ and $\zeta = \left[ \begin{array}{cc} \nu^{\rm T} & w^{\rm T} \end{array} \right]^{\rm T}$. Then,
\begin{equation}
    {z}_{k+1} = \breve{A}_k z_k + \breve{\Gamma}_k \zeta_k + \breve{G}_k u_T,
\label{eqn:A1}
\end{equation}
where
$$ \breve{A}_k =\left[\begin{array}{cc} \bar{A}_k & - \eta_k B K \\
\mathbf{0} & \tilde{A}_k
\end{array} \right],~
\breve{\Gamma}_k =\left[\begin{array}{cc} \mathbf{0} & 
\mathbf{I}  \\
(1-\eta_k) L & \mathbf{I}
\end{array} \right],$$
and 
$$\breve{G}_k = \left[\begin{array}{c} \mathbf{0} \\
\eta_k B \end{array} \right]. $$

Under contractivity conditions (\ref{eqn:def_1a}) and (\ref{eqn:def_1b}), repeat the analysis in Propositions 1 and 2 for (\ref{eqn:A1}). Let $\breve{P}_k = \mathbb{E} [(z_k - \mu_{z,k}) (z_k - \mu_{z,k})^{\rm T}]$, where $\mu_{z,k} = \mathbb{E}[z_k]$. Then,
\begin{equation}
\breve{P}_{k+1}=\breve{A}_k \breve{P}_{k}\breve{A}_k^{\rm T}
+ \breve{\Gamma}_k \left[\begin{array}{cc} \Sigma_\nu & 
\mathbf{0}  \\
\mathbf{0} & \Sigma_w
\end{array} \right] \breve{\Gamma}_k^{\rm T},
\label{eqn:A2}
\end{equation}
and $\breve{P}_k \to \breve{P}_k^{\rm s}$ as $k \to \infty$,
where $\breve{P}_k^{\rm s}$ is the unique $N$-periodic solution to (\ref{eqn:A2}).  Then, the steady-state covariance matrix satisfies
$\breve{P}_{x,k}^{\rm s} = \left[\begin{array}{cc} \mathbf{I} 
& \mathbf{0} \end{array}
\right]  \breve{P}_k^{\rm s} 
\left[\begin{array}{cc} \mathbf{I} & \mathbf{0} \end{array}
\right]^{\rm T}.$
Thus, $x_k$ converges to $x_k^{\rm s}$, a cyclostationary process with the $N$-periodic mean, $\mu_{x,k}^{\rm s}$, and $N$-periodic covariance matrix, $P_{x,k}^{\rm s}$.

As no assumption on the actual probability density functions of $w_k$ and $\nu_k$ is made, we resort to the multivariate Chebyshev's inequality \cite{Cheb1, Cheb2} to treat the chance constraint, which states
\begin{equation}
    \text{Prob} \bigg(( x_k-\mu_{x,k}^{\rm s} )^{\rm T} (P_{x,k}^{\rm s})^{-1} (x_k-\mu_{x,k}^{\rm s}) \leq \alpha_x^2 \bigg) \geq 1-\frac{n_x}{\alpha_x^2},
\label{eqn:cheby1}
\end{equation}
where $n_x$ is the dimension of $x$ and $0 < \alpha_x \leq 1$.

Given $\delta > 0$, choose $\alpha_x = \sqrt{n_x / \delta}$. Suppose that for $k = 0, 1, \dots, N-1$, the following condition is verified:
\begin{equation}
    \mu_{x,k}^{\rm s} \in X \sim \mathcal{E} \left( \mathbf{0},\frac{1}{\alpha_x^2} (P_{x,k}^{\rm s})^{-1} \right),
\label{eqn:cheby2}
\end{equation}
where $\mathcal{E}(\mathbf{0},S)=\{x:~x^{\rm T} S x \leq 1\}$ is an ellipsoidal set and $\sim$ denotes the Pontryagin set difference. Then, in steady-state, the chance constraint (\ref{eqn:cnc_con}) holds.

{\bf Remark 6:} 
The steady-state periodic solution, $\breve{P}_{k}^{\rm s}$, can be computed by solving the conventional discrete-time Lyapunov equation for the evolution of the {\em lifted} system error covariance matrix, i.e., of $\breve{P}^l=\text{diag}\{\breve{P}_{0}^{\rm s},\cdots,\breve{P}_{N-1}^{\rm s}\}$, which is directly obtained from (\ref{eqn:A2}).

\section{Selecting an Optimal Admissible Sequence}
\label{sec:opt}

Since multiple admissible sequences may exist, we can select one by minimizing a cost functional. Consider the blended cost $J$ that penalizes the estimation, the control objective, and the control effort:
\begin{equation}
    J = \frac{1}{N} \sum_{k = 0}^{N-1} \bigg( \text{Tr} \left( R_e P_k^{\rm s} \right) + 
    \text{Tr} \left( R_x P_{x,k}^{\rm s} \right) + r_{\eta} \eta_k \bigg),
\label{eqn:cost1}
\end{equation}
where $R_e = R_e^{\rm T} \succeq 0$, $R_x = R_x^{\rm T} \succeq 0$ and $r_{\eta} \geq 0$ are weights. The $1/N$ factor in the cost functional normalizes for sequence length, ensuring that a given reducible/irreducible sequence pair yield the same cost.



Now, search over sequences of a fixed length, check for admissibility, then find the admissible sequence that yields the lowest cost. If no admissible sequence exists, we then extend the sequence length and search again. The process is summarized by the following algorithm:
\begin{algorithm}
\caption{Sequence search up to length $N$}
\begin{algorithmic}[1]
\STATE Fix sequence length $N \in \mathbb{Z}_{>0}$.

\STATE Form $2^N$ binary integer sequences.

\STATE For each sequence, determine the associated irreducible subsequence.

\STATE Check the conditions of admissibility for the irreducible subsequence using dwell-time conditions (\ref{eqn:ncc2restated}), (\ref{eqn:ncc2d2xrestated}) or directly based on (\ref{eqn:def_1a}), (\ref{eqn:def_1b}).

\STATE If subsequence is admissible, evaluate the cost functional (\ref{eqn:cost1}).

\STATE If no such subsequence satisfies the conditions of admissibility, then increase $N$ and return to Step 1.

\STATE If at least one such subsequence is found to be admissible, then select the sequence that minimizes the objective function $J$.
\end{algorithmic}
\end{algorithm}

The equivalence of admissibility between a sequence and its associated irreducible subsequence is invoked after incrementing $N$ in Step 6; if Step 3 produces an irreducible subsequence that has already been evaluated in a previous iteration, then it can be skipped. 

\section{Relative Motion Control Example}
\label{sec:opt_cost_ex}

We consider a case study of spacecraft three dimensional relative motion control.  The relative motion dynamics are modeled with the linearized Clohessy-Wiltshire \cite{CW} equations,
\begin{equation}
\begin{split}
    \ddot{x}_1 &= 3 \omega^2 x_1  + 2 \omega \dot{x}_2 + \frac{1}{m} u_1, \\
    \ddot{x}_2 &= -2 \omega \dot{x}_1 + \frac{1}{m} u_2,\\
    \ddot{x}_3 &= - \omega^2 x_3 + \frac{1}{m} u_3,
\end{split}
\end{equation}
where $m$ is the chaser vehicle's mass and $\omega$ is the mean motion of the target vehicle's orbit. We form a system of first-order equations with state vector $x_k = [ x_{1,k}, x_{2,k}, x_{3,k}, \dot{x}_{1,k}, \dot{x}_{2,k}, \dot{x}_{3,k} ]^\textrm{T}$ and discretize using a Zero-Order Hold \cite{Toth} with sampling period of $30 ~sec$, chaser vehicle mass of $140 ~kg$ and target vehicle mean motion of $\omega = 0.0010 ~rad/sec$.
We choose $C = [ \mathbf{I}_3 \ \ \mathbf{0}_{3 \times 3} ]$, which corresponds to relative position measurements. The goal in this scenario is to rendezvous the chaser vehicle with the target vehicle, i.e., to bring the chaser's state to the origin.

For this controller and observer, the feedback gain matrix $K$ is computed using LQR by solving the Discrete Time Algebraic Riccati Equation (DARE) and the observer gain matrix $L$ computed by solving the dual DARE, with $Q = \mathbf{I}_6$ and $R = \mathbf{I}_3$ in each case. 


The distributions of the measurement noise, $\nu_k$, and of the state disturbance, $w_k$, are assumed to be Gaussian with covariance matrices $\Sigma_v = 10^{-2} \cdot \mathbf{I}_3$ and $\Sigma_w = 10^{-4} \cdot \mathbf{I}_6$.


The cost function (\ref{eqn:cost1}) has been defined with $r_\eta = 0$, $R_e = \mathbf{I}$ and $R_x = \mathbf{0}$. These choices ensure accurate estimates of the relative position states.

We  now construct a sequence that satisfies the sufficient conditions, and demonstrate that it leads to stable behavior.

For this example, $\rho(\bar{\Omega}_0) = 1.0063$, $\rho(\bar{\Omega}_1) = 0.2016$, $\rho(\tilde{\Omega}_0) = 0.0332$, and $\rho(\tilde{\Omega}_1) = 1.0063$.
Working with the Frobenius norm, $k^* = 1$ and
\[
    c = \frac{\|\bar{\Omega}_1\|_F}{\rho(\bar{\Omega}_1)} = \frac{10.4716}{0.2016} = 51.950.
\]
Begin with the (inadmissible) sequence $s_2 = \{ 0, 1 \}$, with $n_s = 2$, $n_0 = 1$, and $n_1 = 1$. Then, $n_s \log c + n_0 \log \bar{\rho}_0 + n_1 \log \bar{\rho}_1 = 6.3054$ and
$n_s \log c + n_1 \log \tilde{\rho}_0 + n_0 \log \tilde{\rho}_1 = 4.5016$.

\begin{figure}
    \centering
    \includegraphics[scale = 0.25]{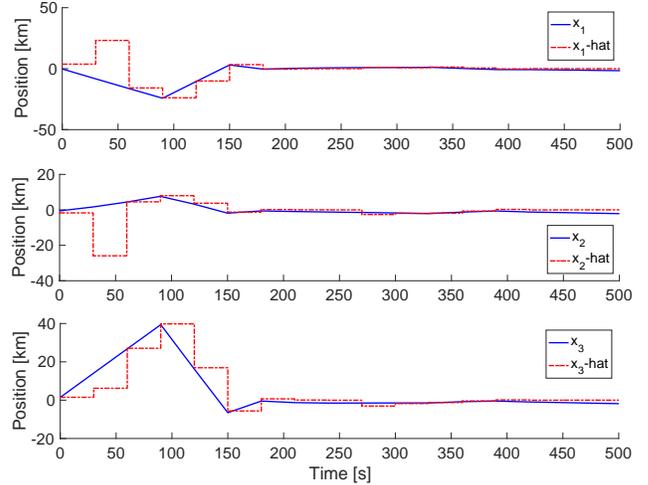}
    \caption{Simulation of the 8-step sequence $s_{8}$ that was constructed to satisfy the sufficient conditions.}
    \label{fig:relmot_8}
\end{figure}

To satisfy the sufficient conditions, increase $n_1$ until (\ref{eqn:ncc2restated}) is satisfied, then increase $n_0$ until (\ref{eqn:ncc2d2xrestated}) is satisfied, then repeat as necessary until both inequalities are simultaneously satisfied. This process yields $n_0 = 3$ and $n_1 = 5$, as then,
$n_s \log c + n_0 \log \bar{\rho}_0 + n_1 \log \bar{\rho}_1 = -0.0879$ and 
$n_s \log c + n_1 \log \tilde{\rho}_0 + n_0 \log \tilde{\rho}_1 = -2.2837$.
Then, we simulate this new sequence, $s_{8} = \{ 0, 0, 0, 1, 1, 1, 1, 1 \}$, that we have constructed, using randomized initial conditions. An example trajectory appears in Figure \ref{fig:relmot_8}.

\begin{figure}
    \centering
    \includegraphics[scale = 0.22]{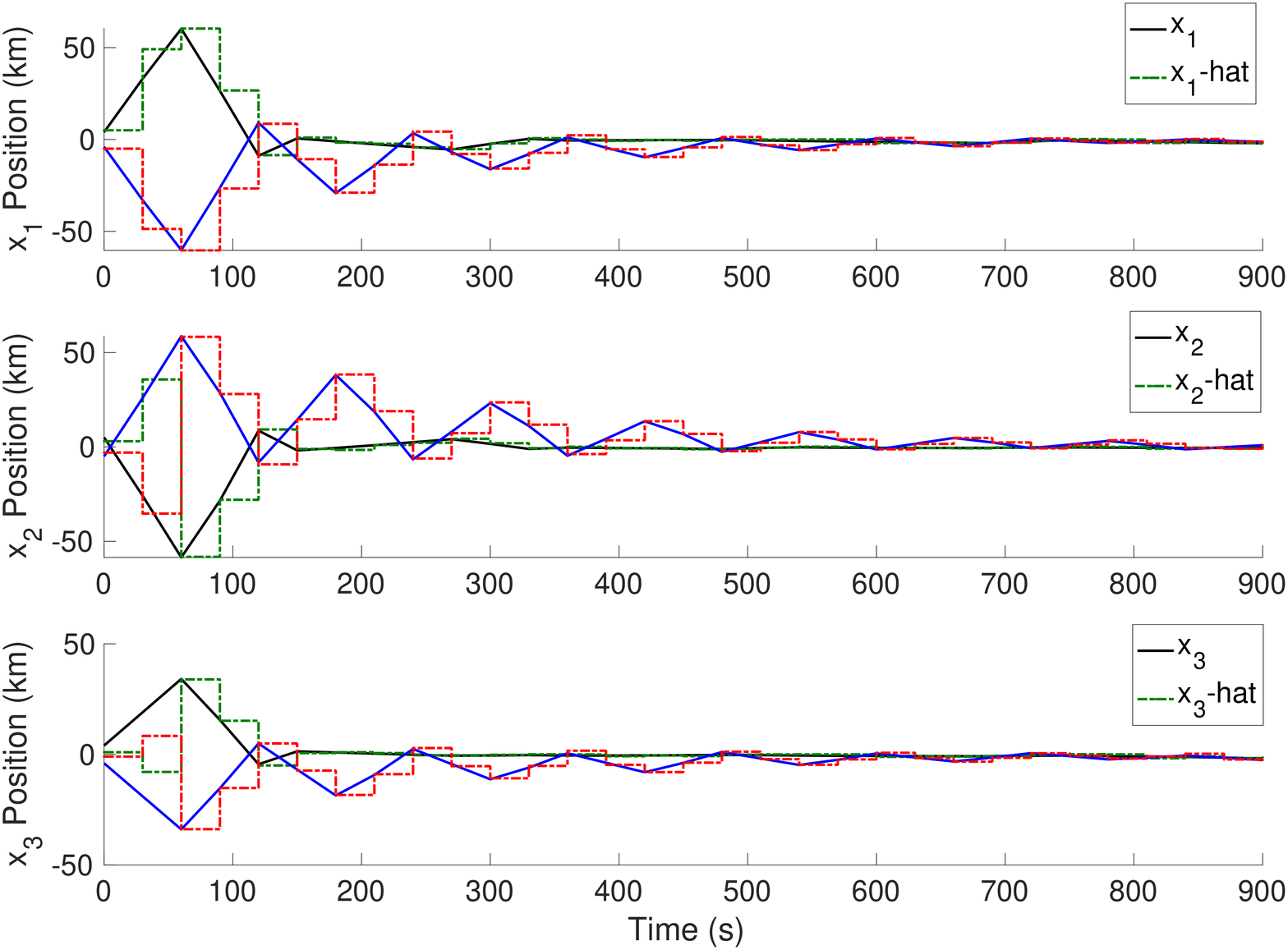}
    \caption{State responses when propagating the length 4 and length 7 optimal sense-control admissible sequences for the relative motion scenario.}
    \label{fig:rel_mot4and7}
\end{figure}

The shortest admissible solution sequence is of length $4$, and the optimum such sequence is $S_4 = \{ 0, 0, 1, 1\}$, with $\bar{q}_A = 0.5879$ in (\ref{eqn:def_1a}) and $\tilde{q}_A = 0.0130$ in (\ref{eqn:def_1b}). When the algorithm treats sequences of length 7, the optimum sequence is $S_7 = \{ 0, 0, 1, 1, 1, 0, 0\}$, with $\bar{q}_A = 0.07594$ and $\tilde{q}_A = 3.796 \times 10^{-5}$. 
Figure \ref{fig:rel_mot4and7} shows the results of propagating these control sequences forward.

In each case, the expected value of each state is successfully driven to the origin. Of note, when the algorithm treats sequences of length 8, the optimum sequence is $\{ 0, 0, 1, 1, 0, 0, 1, 1\}$, a reducible sequence which has $S_4$ as its corresponding irreducible subsequence.

\begin{figure}
    \centering
    \includegraphics[scale = 0.215]{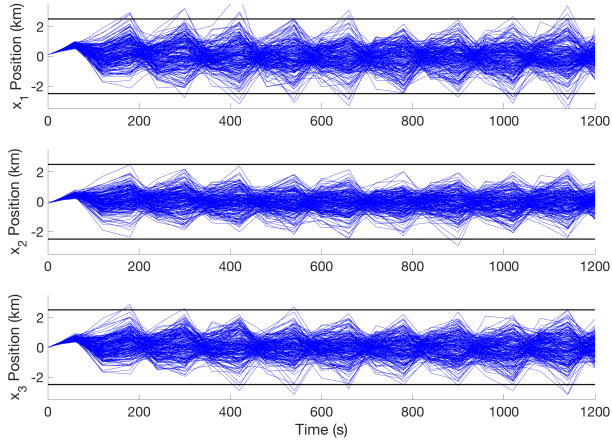}
    \caption{Two hundred simulated trajectories of the relative position of the chaser vehicle under the sequence $S_4$, subject to the box constraint $b = 2.5$ km and the chance constraint P$(X) = 0.95$.}
    \label{fig:200_sims}
\end{figure}

The remaining objective is to select an admissible sequence that also satisfies a specified chance constraint. For the relative motion scenario, suppose that we wish to establish a constraint on the steady-state of the form $X = \{ x : \| x \|_\infty \leq b \}$, so that the chaser spacecraft remains within a box centered at the origin, of side length $2b$, with Prob$(X) \geq 0.95$. Consider again sequence $S_4$; invoking (\ref{eqn:cheby2}), with $\alpha_x = \sqrt{3 / 0.05}$, for each of $P^s_{x,0}, \dots, P^s_{x,3}$ yields spheres with a minimum radius of $b = 2.79$ and a maximum radius of $b = 9.54$. Only the loosest of these constraints is guaranteed by (\ref{eqn:cheby2}), but even the tightest bound remains somewhat conservative. When we set $b = 2.5$, as demonstrated in Figure \ref{fig:200_sims}, the chance constraints are violated in no more than $4 \%$ of trajectories at any given time step over the course of two hundred simulation runs.

\section{Concluding remarks}
\label{sec:concl}

This paper formulated and treated a problem in which simultaneous sensing and actuation are not possible. The solution involved the use of offline-constructed periodic switching sequences between sensing and actuation. When applied online, these sequences had desirable convergence properties. Approaches to simplify the check for admissibility of a sequence have been described based on the notion of reducible sequences and the use of dwell time sufficient conditions. An example of applying the procedure to spacecraft relative motion control has been given in simulations.

\end{document}